\documentclass[11pt]{article}
\usepackage{amsfonts,amssymb, graphicx,a4,bm}
\usepackage[english]{babel}
\usepackage{amsmath,amsthm,epsfig}

\newtheorem{theorem}{Theorem}

\begin{document}

\def\reff#1{(\protect\ref{#1})}

\let\a=\alpha \let\b=\beta \let\ch=\chi \let\d=\delta \let\e=\varepsilon
\let\f=\varphi \let\g=\gamma \let\h=\eta    \let\k=\kappa \let\l=\lambda
\let\m=\mu \let\n=\nu \let\o=\omega    \let\p=\pi \let\ph=\varphi
\let\r=\rho \let\s=\sigma \let\t=\tau \let\th=\vartheta
\let\y=\upsilon \let\x=\xi \let\z=\zeta
\let\D=\Delta \let\F=\Phi \let\G=\Gamma \let\L=\Lambda \let\Th=\Theta
\let\O=\Omega \let\P=\Pi \let\Ps=\Psi \let\Si=\Sigma \let\X=\Xi
\let\Y=\Upsilon

\global\newcount\numsec\global\newcount\numfor
\gdef\profonditastruttura{\dp\strutbox}
\def\senondefinito#1{\expandafter\ifx\csname#1\endcsname\relax}
\def\SIA #1,#2,#3 {\senondefinito{#1#2}
\expandafter\xdef\csname #1#2\endcsname{#3} \else
\write16{???? il simbolo #2 e' gia' stato definito !!!!} \fi}
\def\etichetta(#1){(\veroparagrafo.\veraformula)
\SIA e,#1,(\veroparagrafo.\veraformula)
 \global\advance\numfor by 1
 \write16{ EQ \equ(#1) ha simbolo #1 }}
\def\etichettaa(#1){(A\veroparagrafo.\veraformula)
 \SIA e,#1,(A\veroparagrafo.\veraformula)
 \global\advance\numfor by 1\write16{ EQ \equ(#1) ha simbolo #1 }}
\def\BOZZA{\def\alato(##1){
 {\vtop to \profonditastruttura{\baselineskip
 \profonditastruttura\vss
 \rlap{\kern-\hsize\kern-1.2truecm{$\scriptstyle##1$}}}}}}
\def\alato(#1){}
\def\veroparagrafo{\number\numsec}\def\veraformula{\number\numfor}
\def\Eq(#1){\eqno{\etichetta(#1)\alato(#1)}}
\def\eq(#1){\etichetta(#1)\alato(#1)}
\def\Eqa(#1){\eqno{\etichettaa(#1)\alato(#1)}}
\def\eqa(#1){\etichettaa(#1)\alato(#1)}
\def\equ(#1){\senondefinito{e#1}$\clubsuit$#1\else\csname e#1\endcsname\fi}
\let\EQ=\Eq
\def\0{\emptyset}

\def\pp{{\bm p}}\def\pt{{\tilde{\bm p}}}


\def\\{\noindent}
\let\io=\infty

\def\VU{{\mathbb{V}}}
\def\EE{{\mathbb{E}}}
\def\GI{{\mathbb{G}}}
\def\TT{{\mathbb{T}}}
\def\C{\mathbb{C}}
\def\CC{{\mathcal C}}
\def\II{{\mathcal I}}
\def\LL{{\cal L}}
\def\RR{{\cal R}}
\def\SS{{\cal S}}
\def\NN{{\cal N}}
\def\HH{{\cal H}}
\def\GG{{\cal G}}
\def\PP{{\cal P}}
\def\AA{{\cal A}}
\def\BB{{\cal B}}
\def\FF{{\cal F}}
\def\v{\vskip.1cm}
\def\vv{\vskip.2cm}
\def\gt{{\tilde\g}}
\def\E{{\mathcal E} }
\def\I{{\rm I}}
\def\rfp{R^{*}}
\def\rd{R^{^{_{\rm D}}}}
\def\ffp{\varphi^{*}}
\def\ffpt{\widetilde\varphi^{*}}
\def\fd{\varphi^{^{_{\rm D}}}}
\def\fdt{\widetilde\varphi^{^{_{\rm D}}}}
\def\pfp{\Pi^{*}}
\def\pd{\Pi^{^{_{\rm D}}}}
\def\pbfp{\Pi^{*}}
\def\fbfp{{\bm\varphi}^{*}}
\def\fbd{{\bm\varphi}^{^{_{\rm D}}}}
\def\rfpt{{\widetilde R}^{*}}

\def\tende#1{\vtop{\ialign{##\crcr\rightarrowfill\crcr
              \noalign{\kern-1pt\nointerlineskip}
              \hskip3.pt${\scriptstyle #1}$\hskip3.pt\crcr}}}
\def\otto{{\kern-1.truept\leftarrow\kern-5.truept\to\kern-1.truept}}
\def\arm{{}}
\font\bigfnt=cmbx10 scaled\magstep1

\newcommand{\card}[1]{\left|#1\right|}
\newcommand{\und}[1]{\underline{#1}}
\def\1{\rlap{\mbox{\small\rm 1}}\kern.15em 1}
\def\ind#1{\1_{\{#1\}}}
\def\bydef{:=}
\def\defby{=:}
\def\buildd#1#2{\mathrel{\mathop{\kern 0pt#1}\limits_{#2}}}
\def\card#1{\left|#1\right|}
\def\proof{\noindent{\bf Proof. }}
\def\qed{ \square}
\def\trp{\mathbb{T}}
\def\trt{\mathcal{T}}
\def\Z{\mathbb{Z}}
\def\be{\begin{equation}}
\def\ee{\end{equation}}
\def\bea{\begin{eqnarray}}
\def\eea{\end{eqnarray}}


\title {Improved bounds  on  coloring of graphs}

\author{Sokol Ndreca$^1$,   Aldo Procacci$^2$, Benedetto Scoppola$^3$ \\\\\small{$^1$Dep. Estat\' istica-ICEx, UFMG, CP 702 Belo Horizonte - MG, 30161-970 Brazil}
\\\small{$^2$Dep. Matem\'atica-ICEx, UFMG, CP 702 Belo Horizonte - MG, 30161-970 Brazil}\\
\small{$^3$Dipartimento di Matematica - Universita Tor Vergata di Roma, 00133 Roma, Italy}\\\\\footnotesize{emails:~~~{sokol@est.ufmg.br};~} { aldo@mat.ufmg.br};~~scoppola@mat.uniroma2.it}
\date{}
\maketitle

\def\lc{\lceil}
\def\rc{\rceil}

\numberwithin{equation}{section}

\begin{abstract}
Given a graph  $G$ with maximum degree $\Delta\ge 3$,
we prove that the
acyclic edge chromatic number $a'(G)$ of $G$ is such that $a'(G)\le\lceil 9.62 (\Delta-1)\rceil$. Moreover we prove that:
 $a'(G)\le \lceil 6.42(\Delta-1)\rceil$ if  $G$  has girth $g\ge 5\,$;   $a'(G)\le \lceil5.77 (\Delta-1)\rc$ if
 $G$  has girth $g\ge 7$;  $a'(G)\le \lc4.52(\D-1)\rc$ if $g\ge 53$;
 $a'(G)\le \D+2\,$ if $g\ge  \lceil25.84\D\log\D(1+ 4.1/\log\D)\rceil$.
 We further prove that the  acyclic (vertex) chromatic number $a(G)$ of $G$  is such that
 $a(G)\le \lc 6.59 \Delta^{4/3}+3.3\D\rc$.
We also prove that the  star-chromatic number $\chi_s(G)$ of $G$ is such that $\chi_s(G)\le
\lc4.34\Delta^{3/2}+ 1.5\D\rc$.
We finally prove that the $\b$-frugal
chromatic number $\chi^\b(G)$ of $G$ is such that
$\chi^\b(G)\le \lc\max\{k_1(\b)\D,\; k_2(\b){\D^{1+1/\b}/ (\b!)^{1/\b}}\}\rc$, where $k_1(\b)$ and $k_2(\b)$ are decreasing functions of $\b$ such that
$k_1(\b)\in[4, 6]$ and
$k_2(\b)\in[2,5]$.
 To obtain these results  we use an improved version of the Lov\'asz Local Lemma due to Bissacot, Fern\'andez, Procacci and Scoppola \cite{BFPS}.
\end{abstract}

\numsec=1\numfor=1
\section{Introduction}
Let $G=(V,E)$ be an undirected graph with vertex set $V$ and edge set $E$.  Let $\D$ be the maximum degree
 of $G$ and  $g$ the girth of $G$
(i.e. the length of the shortest cycle in $G$).
A vertex coloring of $G$ is  {\it proper} if no two adjacent vertices receive the same color.
A proper vertex coloring of $G$ is
{\it acyclic} if there are no two-colored cycles in $G$. A proper vertex coloring of $G$ is
a {\it star coloring} if no path of  length 3 is bi-chromatic.
A proper vertex coloring of $G$ is {\it $\b$-frugal} if any vertex has at most
$\b$  members of any color class in its neighborhood.
Similarly, an edge coloring of $G$ is said to be proper if no pair of
incident edges receive the same color. A proper edge coloring of $G$ is said to be
acyclic if there are no two-colored cycles.

The minimum number of colors required such that a graph $G$ has
at least one proper vertex coloring is called {\it chromatic number} of $G$ and will be denoted by $c(G)$.
The minimum number of colors required such that a graph $G$ has
at least one acyclic proper vertex coloring is called {\it acyclic chromatic number} of $G$ and will be denoted by $a(G)$.
The minimum number of colors required for a graph $G$ to have
at least one star vertex coloring is called the {\it star  chromatic number} of $G$ and will be denoted by $\chi_s(G)$.
The minimum number of colors required such that a graph $G$ has
at least one $\b$-frugal proper vertex coloring is called the {\it $\b$-frugal chromatic number} of $G$ and
will be denoted by $\chi^\b(G)$. The minimum number of colors such that a graph $G$ has
at least  one proper edge coloring is called the chromatic edge number  and
will be denoted by $c'(G)$.
The minimum number of colors required for  a graph $G$ to have
at least one acyclic proper edge coloring is called {\it acyclic edge chromatic number} of $G$ and will
be denoted by $a'(G)$.

As far as we know, the best known upper bound for $a(G)$ in graphs with maximum degree $\D$ has been given in \cite{AMR} (see there Proposition 2.2),
where it is proved that
$a(G)\le 50 \D^{4/3}$ for all $\D\ge1$.  However in \cite{AMR} authors  remarked that the constant 50  is not optimal.
The best known upper bound for $a'(G)$
in a graph with maximum degree $\D$ was obtained in \cite{MR} (see there Theorem 2.2)  where it is proved that $a'(G)\le 16\D$ for all $\D\ge 1$.  Recently
such bound has been sensibly improved in \cite{MNS} if one excludes graphs with girth less than  9. Actually  it is proved in \cite{MNS}
 that, if  $g\ge 9$ and $\D\ge 4$, then $a'(G)\le 5.91\D$   and if $g\ge 220$ and $\D\ge 4$ then  $a'(G)\le 4.52\D$ (see there Theorems 1 and 2).
Alon, Sudakov and Zaks have conjectured
in \cite{ASZ} that $a'(G)\le \D+2$ and they proved this conjecture
for graphs with girth $g\ge 2000\D\log\D$ and $\D\ge 3$ (see there Theorem 4).  Also  in this  case authors did not try to optimize  the constant.
Recently,
the conjecture that $a'(G)\le \Delta + 2$  has been confirmed for some more families of graphs.
Namely, complete bipartite graphs \cite{BS},
outerplanar graphs \cite{MNS2},  and graphs with maximum degree four \cite{BS2}.
To our knowledge, the best known upper bound for $\chi_s(G)$ in graphs of maximum degree $\D$ has been
given in \cite{FRR}, where it is proved that
$\chi_s(G)\le 20 \D^{3/2}$ for $\D\ge 1$ (see there Theorem 8.1).
Finally, Hind, Molloy and Reed have proved in \cite{HMR} that
 $\chi^\b(G)\le \max\{(\b+1)\D, e^3{\D^{1+1/\b}/ \b}\}$ for sufficiently large $\D$ (see there Theorem 2). In papers \cite{AMR}, \cite{ASZ}, \cite{FRR}, \cite{HMR}, \cite{MR}, \cite{MNS}
 the proofs rely on the Lov\'asz Local Lemma.

The Lov\'asz Local Lemma, which  is one of the main  tools of the probabilistic method in combinatorics, has been recently related  to the cluster expansion
of the abstract polymer gas, which in turn is a widely used technique in statistical mechanics.
Indeed, during the   last decade, the intersection between statistical mechanics and combinatorics has attracted the attention of several researchers
and has been increasingly investigated. In particular,
the application of cluster expansion methods to  coloring problems in graph theory goes back to 2001 with a  seminal paper by Sokal
\cite{S} relating the anti-ferromagnetic Potts model partition function on  a graph $G$ with the chromatic polynomial  on the same graph.

Concerning specifically the Lov\'asz Local Lemma, its  surprising and close
connection with statistical mechanics was pointed out by
Scott and Sokal \cite{SS} in 2005. Indeed, in \cite{SS}, using
an old theorem by Shearer \cite{Sh}, the authors showed that the conclusions of the Lov\'asz Local Lemma hold for the
dependency graph $G$ with vertex set $X$ and probabilities $\{p_x\}_{x\in X}$ if and only if the independent-set
polynomial for $G$ is non vanishing in the polydisc of radii $\{p_x\}_{x\in X}$. The relation with statistical mechanics
occurs because the independent-set
polynomial of $G$ { is}, modulo a constant, the partition function of the hard core self repulsive lattice gas on $G$, so that
its logarithm is the pressure of such a gas.
From this, Sott and Sokal could conclude that the Lov\'asz Local Lemma is  a different way to rephrase
the Dobrushin condition \cite{Dob} for  the convergence of  the pressure of the hard core lattice gas on $G$.

In 2007,  Fern\'andez and Procacci  \cite{FP} provided a new   criterion for  the convergence of  the pressure of the lattice gas on a graph  $G$,
and showed  that this  new criterion is always more effective
than the Dobrushin's criterion. Later, the same authors used their criterion in \cite{FP2} to improve Sokal's results on zero-free regions of chromatic polynomial of \cite{S}.

Very recently,  Bissacot {\it et al.} \cite{BFPS} used the Fern\'andez-Procacci criterion \cite{FP} and the results in \cite{SS}  to
improve the Lov\'asz Local Lemma. This new version of the  Lov\'asz Local Lemma has already been  used to improve
an old result on Latin-transversal \cite{BFPS} and
to obtain some new results  about colorings   of the edges of the complete graph $K_n$ \cite{BKP}.
Finally, it is worth mentioning  that in a very recent paper,   Pegden  \cite{Pe} has shown that the new  Lemma  of \cite{BFPS} also holds in the  Moser and Tardos's  algorithmic contest \cite{MT}.

The present paper thus aims  at informing  the  combinatorics community
that many classical bounds  obtained via the Lov\'asz Local Lemma  can be improved using the new version presented in  \cite{BFPS}. In an effort to  convince the reader about that,  we  focus
our attention on  the six  graph coloring problems described above, showing  how it is possible to obtain,  in a quite straightforward way,
improvements  on the bounds  given  in \cite{AMR}, \cite{ASZ}, \cite{FRR}, \cite{HMR}, \cite{MR} and \cite{MNS}, by simply
using the new lemma of \cite{BFPS} in place of the Lov\'asz Local Lemma.

The rest of the paper is organized as follows. In Section 2 we recall the  Lov\'asz Local Lemma (Theorem 1), we present  the new lemma \cite{BFPS} (Theorem 2) and we state our results on graph colorings
(Theorem 3). In Section 3 we prove Theorem 3.

\numsec=2\numfor=1
\section{The Lov\'asz Local Lemma, the new lemma, and results}
We first state the  Lov\'asz Local Lemma (LLL)
and immediately after the new lemma of \cite{BFPS} with the intent of
clarifying how the improvement works.

To state these lemmas we need some definitions.
Hereafter, if  $U$ is a finite set,  $|U|$ denotes its cardinality.
Let $X$ be a finite set.
Let $\{A_x\}_{x\in X}$ be a family of events on some probability space,
each of which having probability $\mathbb {\rm Prob}(A_x)=p_x$ to occur.
A graph $H$ with vertex set $V(H)=X$
is a {\it dependency graph} for the family of events
$\{A_x\}_{x\in X}$  if, for each $x\in X$,
$A_x$ is independent of all the events in the $\sigma$-algebra
generated by $\{ A_y: y \in
{X}\backslash \Gamma^{*}_H(x) \}$, where $\Gamma_H(x)$ denotes the set of
vertices of $H$ adjacent to $x$ and
$\Gamma^{*}_H(x)=\Gamma_H(x)\cup\{x\}$.

 Denoting $\bar A_x$  the complement event of $A_x$,  let $\bigcap_{x\in X}\bar A_x$ be the event such that none of the events
$\{A_x\}_{x\in {X}}$  occurs.

\begin{theorem}[{\bf Lov\'asz Local Lemma}]
\label{LLL} Suppose that $H$ is a dependency
graph for the family of events $\{A_x\}_{x\in X}$ each one  with probability
${\rm Prob}(A_x)=p_x$ and there exist $\{\m_x\}_{x\in X}$ real
numbers in $[0,+\infty)$ such that, for each $x\in X$,
$$
p_x \;\leq\; \; {\mu_x \over \varphi_x(\bm \mu)}\Eq(dobc)
$$
with\
$$
\varphi_x(\bm \mu)=1+\sum_{R\subseteq  \Gamma^*_H(x)}
\prod_{x\in R}\mu_x\Eq(dobr)
$$
 Then
$$
{\rm{Prob}}(\bigcap_{x\in X}\bar A_x)\,> \,0\Eq(a)
$$
\end{theorem}
\\{\bf Remark}. In the literature the LLL is usually written in terms of variables $r_x= {\m_x/( 1+\m_x)}\in [0,1)$,
so that
\equ(dobc)
becomes  $p_x\leq{r_x \prod_{y\in \G_H(x)}( 1-r_y)}$ (see e.g. Lemma 5.1.1 p. 68 in \cite{AS}).
However, the  formulation above, although completely equivalent to the usual one,
shows in a clear way the difference and the consequent  improvement contained in the following theorem.

\begin{theorem}[{\bf\cite{BFPS}}]\label{LLLn} Suppose that $H$ is a dependency
graph for the family of events $\{A_x\}_{x\in X}$ each one with probability
${\rm Prob}(A_x)=p_x$ and there exist $\{\mu_x\}_{x \in X}$ real
numbers in $[0,+\infty)$ such that, for each $x\in X$,
$$
p_x\;\leq\;  {\mu_x \over \ffp_x(\bm \mu)}
\Eq(eq:r4)
$$
with\
$$
\ffp_x(\bm \mu)=1+\sum_{{R\subseteq  \Gamma^*_H(x)}\atop R\ {\rm indep\ in}\ H}
\prod_{x\in R}\mu_x
\Eq(ffp)
$$
 Then
$$
{\rm{Prob}}(\bigcap_{x\in X}\bar A_x)\,> \,0
\Eq(eq:r3b)
$$
\end{theorem}

\\{\bf Remark}. The only difference between the LLL (as stated in Theorem \ref{LLL}) and  Theorem \ref{LLLn} above
is that in Theorem \ref{LLL} the sum of the right hand side of \equ(dobr) is over all the subsets of $\Gamma^*_H(x)$  while
in Theorem \ref{LLLn}  the same sum is now only over  the {\it independent} subsets of  $\Gamma^*_H(x)$.
This yields  $\varphi^*_x(\bm \mu) \le \varphi_x(\bm\mu)$ so that condition \equ(eq:r4) in Theorem \ref{LLLn}
is always less restrictive than  \equ(dobc) in Theorem \ref{LLL}. Moreover, noting that
 $\varphi_x(\bm \mu)=(1+\m_x)\prod_{y\in \G_H(x)}( 1+\m_y)$, it is clear that condition \equ(dobc) of Theorem \ref{LLL}
 does not depend on the graph structure  of $\G_H(x)$  (i.e. on the subgraph of $H$ induced by $\G_H(x)$) but only on its cardinality.
 For example, condition \equ(dobc) is the same, either if
$\G_H(x)$ is an independent set, or
$\G_H(x)$ is a clique.
 In contrast, condition
\equ(eq:r4) in Theorem \ref{LLLn}  does depend on the graph structure  of $\G_H^*(x)$ (and hence of $\G_H(x)$). Consequently, the
improvement brought by Theorem \ref{LLLn} is maximal when the set of vertices $\G_H(x)$ which are neighbors of $x$
form
a clique, and it is nearly null when  vertices of $\G_H(x)$ form an independent set in the dependency graph (e.g. like in bipartite graphs).  In view of this,
in the next section we will frequently the following inequality. Let $x$ be
a vertex of the dependency graph $H$ for the events
$\{A_x\}_{x\in X}$, and suppose that   $\G^*_H(x)$ is the union  (not necessarily disjoint) of  $c_1,\dots ,c_k$ cliques, then, by definition \equ(eq:r4),

$$
\ffp_x(\bm \mu)\le 1+ \sum_{s=1}^k~~\sum_{1\le\, i_1\,<\cdots<\,i_s\le \,k}~~\sum_{x_1\in c_{i_1}}\cdots
\sum_{x_s\in c_{i_s}}\m_{x_1}\dots \m_{x_s}
$$
$$
= \prod_{i=1}^k\Big[1+ \sum_{y\in c_i}\m_y\Big]~~~~~~~~~~~~~~~~~~~~~~~~~~~~~~~\Eq(key)
$$

\vskip.2cm
We conclude the section by stating   the results contained in the present paper, which  can be summarized  by  the following theorem.
\begin{theorem}\label{tot}
If $G$ is a graph with maximum degree $\D\ge 3$ and girth $g$,  then

\begin{itemize}
\item[(a)] $a'(G)\le \lceil 9.62 (\Delta-1)\rceil$.

\item[(b)]\label{gir}
If  $g\ge 5$,  then
 $a'(G)\le \lceil 6.42(\Delta-1)\rceil$. If $g\ge 7$, then $a'(G)\le \lceil5.77 (\Delta-1)\rceil$.
If $g\ge 53$, then $a'(G)\le \lceil 4.52(\D-1)\rceil$.

\item[(c)]\label{d+2}
If
$g\ge  \lceil25.84\D\log\D\big(1+ {4.1\over \log\D}\big)\rceil$,
then $a'(G)\le \D+2$.

\item[(d)]\label{large}
$a(G)\le \lceil 6.59\Delta^{4/3}+ 3.3\D\rceil$.

\item[(e)]\label{star}
$\chi_s(G)\le \lceil4.34\Delta^{3/2}+ 1.5\D\rceil$.

\item[(f)]
For any $\b\ge 1$,  $\chi^\b(G)\le \lceil\max\{k_1(\b)\D, \; k_2(\b){\D^{1+{1\over \b}}\over (\b!)^{1/\b}}\}\rceil$, where $k_1(\b)$ and $k_2(\b)$ are decreasing functions of $\b$ such that
$k_1(\b)\in[4, 5.27]$
and
$k_2(\b)\in[2,4.92]$.
\end{itemize}
\end{theorem}
\\{\bf Remark}. Note that, differently from Theorem 2 in \cite{HMR}, in item (f) it is not required  to take $\D$ sufficiently large. As for items (a)\,-(e), to prove item (f) we only need $\D\ge 3$. We also stress that we did not  attempt to optimize the non leading terms in $\D$ in  bounds (c)\,-(e).
\vv

\numsec=3\numfor=1
\section{Proof of Theorem \ref{tot} }
Hereafter $G=(V,E)$ will denote an undirected graph with vertex
set $V$, edge set $E$, maximum degree $\D\ge 3$ and girth $g$.

\subsection{Proof of item (a): acyclic edge chromatic number of $G$}
\def\N{{\mathbb N}}
Let  $K$ be the set whose elements are the pairs $\{e,e'\}\subset E$ such that $e, e'$ are
incident in a common vertex. Let, for  $k\ge 2$, $C_{2k}(G)$
be the set of all cycles in $G$ of length $2k$.  Finally, let $X= K\cup(\,\bigcup_{k\ge 2}C_{2k})$.
We regard cycles $c_{2k}\in C_{2k}$ as sets of edges, so that the elements of $X$ are (some of) the subsets of $E$.
For each edge $e\in E$, choose a color  independently  and uniformly  among $N$ possible colors
such that $N\ge c(\D-1) $ ($c$ is a constant to be determined later). Consider now the following unfavorable events.

\begin{itemize}
\item[I.] For $\{e,e'\}\in K\mathbb{}$, let $A_{\{e,e'\}}$ be the event that  the  edges $e$ and $e'$ have the same color.

\item[II.] For   $c_{2k}\in C_{2k}$ ($k\ge 2$), let $A_{c_{2k}}$ be the event that the cycle $c_{2k}$ is (properly) bichromatic.
\end{itemize}

If condition \equ(eq:r4) of Theorem \ref{LLLn} holds,  there is a non zero probability
that none of the events of type I or II occurs, and hence there exists a proper edge coloring  of $G$  with
no two-colored cycles. To check condition \equ(eq:r4) we first observe that, for each $\{e,e'\}\in K$,  the probability
of the event $A_{\{e,e'\}}$ is
$$
{\rm Prob}(A_{\{e,e'\}}) = {1\over  N}
$$
while, for any $k\ge 2$ and $c_{2k}\in C_{2k}$
$$
{\rm Prob}(A_{c_{2k}})\le {1\over  N^{2k-2}}
$$
Secondly, we have to find a graph with vertex set $X$ which is a dependency graph for the events
$\{A_{x}\}_{x\in X}$.
Since we are choosing a color at random for each edge independently, we have clearly that the event $A_{\{e,e'\}}$ is independent
of any other event $A_{\{f,f'\}}$
such that $\{e,e'\}\cap \{f,f'\}=\emptyset$ and  of all events  $A_{c_{2k}}$ with  $k\ge 2$ such that $\{e,e'\}\cap c_{2k}=\emptyset$.
Analogously, for any $m\ge 2$,
the event $A_{c_{2m}}$ is independent of all events  $A_{\{f,f'\}}$ with  $\{f, f'\}\in K$ and all events $A_{\tilde c_{2k}}$ with $k\ge 2$ such that
$c_{2m}\cap \{f,f'\}=\emptyset$ and $c_{2m}\cap \tilde c_{2k}=\emptyset$. So let $H=(X,F)$ be the
graph with vertex set $X$ and edge set $F$ such that
the pair $\{x,x'\}\in F$ if and only if  $x\cap x'\neq\emptyset$. By construction, $H$ is a dependency graph for
the  events $\{A_{x}\}_{x\in X}$. Now observe that
\begin{itemize}
\item
each edge $e$ is contained in at most $2(\D-1)$ pairs $\{f,f'\}\in K$
\item
each edge $e$ is contained in at most $(\D-1)^{2k-2} $ cycles  $c_{2k}\in C_{2k}$,
for any $k\ge2$
\end{itemize}

Hence

\begin{itemize}
\item[{[a]}]
for each  vertex $x=\{e_1,e_2\}\in K$ of $H$,
$\G_H^*(x)$ is the union of
two sets $\G_1^*(x)$ and
$\G_2^*(x)$ such that, for $i=1,2$
$$
|\G_i^*(x)|\le 2(\D-1) + \sum_{s\ge 2}(\D-1)^{2s-2}
$$
and every
element $z\in \G_i^*(x)$ contains $e_i$ ($i=1,2$), so that
the subgraph of $H$  induced by $\G_i^*(x)$ is a clique.
\item[{[b]}]
For $k\ge 2$ and for each  vertex $y=c_{2k}=\{e_1,\dots,e_{2k}\}\in C_{2k}$ of $H$,
$\G_H^*(y)$
is the union of $2k$ sets $\G_1^*(y),\dots,\G_{2k}^*(y)$ such that,
for any $j=1,2,\dots,  2k$,
$$
|\G_j^*(y)|\le 2(\D-1) + \sum_{s\ge 2}(\D-1)^{2s-2}
$$
and every element $z\in \G_{j}^*(y)$ contains $e_j$.
Hence, for any $j=1,2,\dots,  2k$,  the subgraph of $H$  induced by $\G_j^*(y)$ is a clique.
\end{itemize}
Let us now choose nonnegative  numbers $\{\m_z\}_{z\in X}$ such that: for any $x\in K$, $\m_{x}=\m_1$;  for
each $y\in C_{2k}$,
$\m_{y}=\m_k$. Then, recalling definition \equ(ffp) and inequality \equ(key), an easy calculation shows that under conditions [a] and [b] we have

$$
\ffp_{x}(\bm \mu)\;\le\;\Big[1+ 2(\D-1)\m_1 + \sum_{s\ge 2}(\D-1)^{2s-2}\m_s\Big]^2
$$
and
$$
\ffp_{y}(\bm \mu)\;\le\;\Big[1+ 2(\D-1)\m_1 + \sum_{s\ge 2}(\D-1)^{2s-2}\m_s\Big]^{2k}
$$
and hence conditions \equ(eq:r4) become
$$
{1\over N}~\le ~{\m_{1}\over \Big[1+ 2(\D-1)\m_1 + \sum_{s\ge 2}(\D-1)^{2s-2}\m_s\Big]^2}
$$
$$
 {1\over N^{2k-2}}~\le~ {\m_{k}\over \Big[1+ 2(\D-1)\m_1 + \sum_{s\ge 2}(\D-1)^{2s-2}\m_s\Big]^{2k}}
$$
Now choosing $\m_1=\m={\a\over \D-1}$ with $0<\a<1$ and $\m_k=\m^{2k-2}$, and recalling that $N\ge c(\D-1)$,
the conditions above are satisfied if
$$
{1\over c}\le {\a\over (1+ 2\a + \sum_{s\ge 2}\a^{2s-2})^2}
$$
$$
{1\over c}\le {\a\over (1+ 2\a + \sum_{s\ge 2}\a^{2s-2})^{2k/2k-2}}
$$
Since $k\ge 2$,  the first inequality implies the second.  Therefore
the condition which guarantees that none of the bad events $\{A_x\}_{x\in X}$ occurs is
$$
{1\over c}\le {\a\over (1+ 2\a + \sum_{s\ge 2}\a^{2s-2})^2}
$$
i.e.,
$$
{ c}~\ge~ \a^{-1}{\Big[1+ 2\a + {\a^2\over 1-\a^2}\Big]^2}
\Eq(fint)
$$
The function on the right hand side of \equ(fint) can be minimized in the interval $\a\in (0,1)$
and a straightforward calculation  gives that \equ(fint) is satisfied if $c\ge 9.6130002$. Hence every graph $G$
with maximum degree $\D$ such that edges are colored using a number of colors $N$
greater or equal than $ 9.62(\D-1)$ admits an acyclic proper coloring.
$\Box$
\vskip.2cm
\noindent{\bf Remark}. As observed in \cite{MNS}, using the Lov\'asz Local Lemma  one can obtain an upper bound for  the edge chromatic number $c'(G)$ of a graph $G$
 at best $c'(G)\le \lceil 4e\D\rceil $ and  for any $G$
we have clearly that  $a'(G)\ge c'(G)$. We leave to the reader to check that one can obtain $c'(G)\le 4(\D-1)$ using Theorem \ref{LLLn} in place of
Theorem \ref{LLL} and proceeding similarly to the scheme illustrated in the proof of item (a) above.

\subsection{Prof of item (b): acyclic edge chromatic number of $G$ when $g\ge 5$ }

We follow here the strategy described in \cite{MNS}.
Namely we will first consider the following
problem.  Let
 $\eta\ge 2$ be an integer. We want to know the minimum colors  needed   to find a coloring $\CC$ of
  the edges of $G$ such that

\begin{itemize}

\item[1.] In any vertex $v$ of $G$ the number of edges incident to $v$ having the same color is at most $\eta$

\item[2.] There is no properly bichromatic cycle in $G$

\item[3.] There is no monochromatic cycle in $G$

\end{itemize}

Suppose that we are able to prove that, for some $N\in \mathbb{N}$, we find a coloring $\CC$  which  satisfies 1, 2, 3,
using $N$ colors. Then it is also possible to find a coloring $\CC'$ using $N'=\eta N$ colors which is proper and satisfies 2
(i.e. $\CC'$ is an acyclic proper coloring).
Indeed,
just observe that in the coloring $\CC$ the sets of edges with the same color are forests with maximum degree $\eta$
and one needs $\eta$ colors to proper color a forest with maximum degree $\eta$. So if one recolors
each color $c_i$ ($i=1,2,\dots N$) in the coloring
$\CC$ using $c_i^1, c_i^2,\dots, c_i^\eta$  distinct colors in such a way that monochromatic forests disappear,
then one gets a new coloring $\CC'$ in which $N'=\eta N$ colors are used and by construction $\CC'$
is proper and satisfies 2.

Now, we use Theorem 2 to show that if  $N\ge c(\D-1)$ (where $c$ is a constant to be determined),
then the coloring $\CC$ satisfying properties 1-3 exists and hence, in view of the above argument,
if $N'\ge c'(\D-1)$, with ${\bar c}=\eta c$,
then there is an acyclic edge coloring $\CC'$ on $G$ using $N'$ colors.

As we did in the previous subsection, let us
choose for each edge $e\in E$ independently a color at random among $N\ge c(\D-1)$ possible colors.
Let now $K_\eta$ be the set whose elements are sets of edges  $\k_\eta=\{e_1,e_2,\dots,e_{\eta+1}\}\subset E$
all incident to a common vertex.
Let $C_{m}$
($m\ge 3$)  be the set whose elements  are  all cycles $c_m$ in $G$ of length $m$.
Finally, let $X= \{\bigcup_{m\ge 3}C_m\}\cup K_\eta$. We regard cycles  as subsets of edges.
So again the elements of $X$ are (some of) the subsets of $E$.
We now consider the following unfavorable events

\begin{itemize}
\item[I.] For $\k_\eta=\{e_1,\dots,e_{\eta+1}\}\in K_\eta$, let $A_{\k_\eta}$ be the event that
 all edges $e_1,\dots,e_{\eta+1}$ have the same color.

\item[II.] For   $c_{2k}\in C_{2k}$, let $A_{c_{2k}}$ be the event that the cycle $c_{2k}$ is either properly bichromatic
or monochromatic.

\item[III.] For   $c_{2l+1}\in C_{2l+1}$, let $ A_{c_{2l+1}}$ be the event that the cycle $c_{2l+1}$ is monochromatic.
\end{itemize}

Theorem 2 gives
a condition  which guarantees that the probability that none of the events of type I or II or III occurs
is strictly positive and hence
 the existence of a coloring $\CC$ of $G$ with properties 1, 2 and 3 above.

Observe that, for $\k_\eta\in K_\eta$,  the probability
of the event $A_{\k_\eta}$ is
$$
{\rm Prob}(A_{\k_\eta})= {1\over N^\eta}
$$
while, for any $k\ge \lceil g/2\rceil$ and $c_{2k}\in C_{2k}$
$$
{\rm Prob}(A_{c_{2k}})= {1\over N^{2k-2}} \Eq(exp)
$$
and, for any $l\ge \lfloor {\,g/ 2}\rfloor$ and $c_{2l+1}\in C_{2l+1}$
$$
{\rm Prob}(A_{c_{2l+1}})= {1\over N^{2l}}
$$
To prove \equ(exp) just observe that the total number of ways of coloring  an even cycle $c_{2k}$
using $N$ colors is $N^{2k}$ while the number of ways of coloring an even cycle $c_{2k}$
using $N$ colors so that the cycle is either monochromatic or proper bichromatic
is $N+N(N-1)=N^2$, where  $N$ is the number
of different monochromatic ways of coloring the  cycle $c_{2k}$ and $N(N-1)$  is the number of
different proper bichromatic ways of
coloring the cycle $c_{2k}$. So ${\rm Prob}(A_{c_{2k}})= N^2/N^{2k}=1/N^{2k-2}$.

Now  we have to find a graph with vertex set $X$ which is a dependency graph for the events
$\{A_{x}\}_{x\in X}$.
Since we are choosing a color at random for each edge independently, we have once again
that the event $A_{x}$ is independent of all other events $A_{x'}$
such that  $x\cap x'=\emptyset$.  So the graph  $H=(X,F)$, with vertex
set $X$ and edge set $F$ such that
the pair $\{x,x'\}\in F$ if and only if $x\cap x'\neq\emptyset$, is   a dependency graph for the events $\{A_{x}\}_{x\in X}$.
Now observe that
\begin{itemize}
\item
each edge $e$ is contained in at most $2{\D-1\choose \eta}\le 2{(\D-1)^{\eta}\over \eta!}$
distinct sets $\k_\eta=\{e_1,\dots,e_{\eta+1}\}\in K_\eta$;
\item
each edge $e$ is contained in at most $(\D-1)^{m-2}$ cycles  $c_{m}\in C_m$ ($m\ge 3$).
\end{itemize}
Hence we have the following.

\begin{itemize}
\item[{[a]}]
For each  vertex $x=\k_\eta=\{e_1,e_2,\dots,e_{\eta+1}\}\in K_\eta$ of $H$,
$\G_H^*(x)$ is the union of $\eta+1$ sets $\G_i^*(x)$ ($i=1,\dots ,\eta+1\}$)
such that
$$
|\G_i^*(x)|\le 2{(\D-1)^{\eta}\over \eta!} +
\sum_{k\ge \lceil g/2\rceil}(\D-1)^{2k-2}+ \sum_{l\ge \lfloor g/2\rfloor}(\D-1)^{2l-1}
$$
and every element of $\G_i^*(x)$ contains $e_i$.
Hence  the subgraph of $H$  induced by $\G_i^*(x)$ is a clique for $i=1,\dots ,\eta+1$.
\item[{[b]}]
For $m\ge 3$ and for each  vertex $y=c_m=\{e_1,\dots,e_m\}\in C_m$ of $H$, we have that
$\G_H^*(y)$ is the union of $m$ sets $\G_1^*(y),\dots,\G_{m}^*(y)$ such that,
for $j=1,\dots,  m$,
$$
|\G_j^*(y)|\le 2{(\D-1)^{\eta}\over \eta!} +
\sum_{k\ge \lceil g/2\rceil}(\D-1)^{2k-2}+ \sum_{l\ge \lfloor g/2\rfloor}(\D-1)^{2l-1}
$$
and every element of $\G_{j}^*(y)$ contains $e_j$ so that the subgraph of $H$  induced by $\G_j^*(y)$ is a clique
for all $j=1,2,\dots, m$.
\end{itemize}
Let us now choose nonnegative  numbers $\{\m_z\}_{z\in X}$ such that: for any $x\in K_\eta$, $\m_{x}=\m_1$;  for
each $y\in C_m$,
$\m_{y}=\m_m$. Then using once again inequality \equ(key) one gets, under conditions [a] and [b], that

$$
\ffp_{x}\bm \mu)\le\Big[1+ 2{(\D-1)^{\eta}\over \eta!}\m_{1} + \sum_{k\ge \lceil g/2\rceil}(\D-1)^{2k-2}\m_{2k}+ \sum_{l\ge \lfloor g/2\rfloor}\D^{2l-1}\m_{2l+1}\Big]^{\eta+1}
$$
and
$$
\ffp_{y}(\bm \mu)\le\Big[1+ 2{(\D-1)^{\eta}\over \eta!}\m_{1} + \sum_{k\ge \lceil g/2\rceil}(\D-1)^{2k-2}\m_{2k}+ \sum_{l\ge \lfloor g/2\rfloor}(\D-1)^{2l-1}\m_{2l+1}\Big]^{m}
$$
and hence condition \equ(eq:r4) of Theorem \ref{LLLn} becomes
$$
{1\over N^{\eta}}~\le ~{\m_{1}\over \Big[1+ 2{(\D-1)^{\eta}\over \eta!}\m_{1} + \sum\limits_{k\ge \lceil g/2\rceil}(\D-1)^{2k-2}\m_{2k}+ \sum\limits_{l\ge \lfloor g/2\rfloor}(\D-1)^{2l-1}\m_{2l+1}\Big]^{\eta+1}}
$$
$$
 {1\over N^{2k-2}}~\le~ {\m_{2k}\over \Big[1+ 2{(\D-1)^{\eta}\over \eta!}\m_{1} + \sum\limits_{k\ge \lceil g/2\rceil}(\D-1)^{2k-2}\m_{2k}+ \sum\limits_{l\ge \lfloor g/2\rfloor}(\D-1)^{2l-1}\m_{2l+1}\Big]^{2k}}
$$
$$
 {1\over N^{2l}}~\le ~{\m_{2l+1}\over \Big[1+ 2{(\D-1)^{\eta}\over \eta!}\m_{1} + \sum\limits_{k\ge \lceil g/2\rceil}(\D-1)^{2k-2}\m_{2k}+ \sum\limits_{l\ge \lfloor g/2\rfloor}(\D-1)^{2l-1}\m_{2l+1}\Big]^{2l+1}}
$$
Now choose $\m_1=\m^\eta$, $\m_{2k}=\m^{2k-2}$,  $\m_{2l+1}=\m^{2l}$ and $\m={\a\over \D-1}$  with $\a\in (0,1)$. Then,
recalling that $N\ge c(\D-1)$, the conditions above are satisfied if

$$
{1\over c}\le {\a\over \Big[1+ 2{\a^{\eta}\over \eta!} + \sum\limits_{k\ge \lceil g/2\rceil}\a^{2k-2}+{1\over (\D-1)} \sum\limits_{l\ge \lfloor g/2\rfloor}\a^{2l}\Big]^{(\eta+1)/\eta}}
$$
$$
{1\over c}\le {\a\over \Big[1+ 2{\a^{\eta}\over \eta!} + \sum\limits_{k\ge \lceil g/2\rceil}\a^{2k-2}+{1\over (\D-1)} \sum\limits_{l\ge \lfloor g/2\rfloor}\a^{2l}\Big]^{2k/(2k-2)}}
$$
$$
{1\over c}\le {\a\over \Big[1+ 2{\a^{\eta}\over \eta!} + \sum\limits_{k\ge \lceil g/2\rceil}\a^{2k-2}+{1\over (\D-1)} \sum\limits_{l\ge \lfloor g/2\rfloor}\a^{2l}\Big]^{(2l+1)/2l}}
$$

\v
\\If    $g\ge 5$, then $k\ge 3$ and $l\ge 2$ and moreover, if   $\eta\ge2$, the three inequalities are satisfied if

$$
{1\over c}\le {\a\over \Big[1+ 2{\a^{\eta}\over \eta!} +{1\over 1-\a^2}\left({ {\a^{2\lceil g/2\rceil-2}}
+ {1\over \D-1}{\a^{2\lfloor g/2\rfloor}}}\right)\Big]^{(\eta+1)/\eta}}
$$
\v
\\Hence, recalling that ${\bar c}=\eta c$,  observing that   ${ {\a^{2\lceil g/2\rceil-2}}
+ {1\over \D-1}{\a^{2\lfloor g/2\rfloor}}}\le {\D\over \D-1}\a^{2\lceil g/2\rceil-2}$ for all $\D\ge 3$ and all $g\ge 3$, and optimizing with respect to $\a\in (0,1)$, we get
$$
{\bar c}\ge {\eta} \min_{\a\in (0,1)} \a^{-1}\Big[1+ 2{\a^{\eta}\over \eta!} + {\D\over \D-1} {\a^{2\lceil g/2\rceil-2}\over (1-\a^2)}\Big]^{(\eta+1)/\eta}\Eq(abo)
$$
If $g\ge 5$,  $\eta=2$, and $\D\ge 3$ a rough calculation (bounding $\D\over \D-1$ by 3/2 for all $\D\ge  3$)  gives ${\bar c}\ge 6.42$ and hence
$a'(G)\le \lceil6.42(\D-1)\rceil$;
If  $g\ge 7$ and $\eta=2$, we get   ${\bar c}\ge 5.77$ and hence
$a'(G)\le \lceil 5.77(\D-1)\rceil$. Finally, if  $g\ge 53$ and $\eta=3$, we get   ${\bar c}\ge 4.52$ and hence
$a'(G)\le \lceil 4.52\D\rceil$. Note that $a'(G)/\D\le 4.52$ is obtained in \cite{MNS} for
 $g\ge 220$.
$\Box$
\vv
\\{\bf Remark}. Observe that, for fixed $g,\,\eta$ the quantity ${\bar c}$ defined in \equ(abo) slightly decreases as $\D\to \infty$ and one can check that
$\lim_{\D\to\infty} {\bar c}(5,2,\D)\le 6.159\dots$, $\lim_{\D\to\infty}{\bar c}(7,2,\D)=5.654\dots$, $\lim_{\D\to\infty} {\bar c}(53,2,\D)=4.511\dots$, yielding a slight
improvement of $a'(G)$ in all three cases considered as $\D\to\infty$.

\subsection{Proof of  item (c): a class of graphs with maximum degree $\D$ and acyclic edge chromatic number $\le \Delta +2$}
 We follow  \cite{ASZ} using Theorem \ref{LLLn} instead of the Lov\'asz Local Lemma.
By Vizing' Theorem \cite{V}, there exists a proper coloring $\mathcal C$  of the edges of $G$ with $\D+1$
colors. An even cycle is called properly half-monochromatic with respect to the coloring $\mathcal C$,  if  one of its halves (a set of alternate edges)
is monochromatic while the other half is not. Note that a  properly  half-monochromatic cycle is never properly  bichromatic by the coloring $\mathcal C$. Observe also that a cycle
of odd length can  never  be properly bichromatic  by the coloring $\mathcal C$.

Let   $K$ be the set whose elements are all pairs   of adjacent edges  $\{e,e'\}$ of $G$.
Let $B_{2m}$ ($m\ge 2$)  be the set whose elements  are  all cycles $b_{2m}$ in $G$ of length $2m$ which are properly bichromatic by the coloring $\mathcal C$.
Let $H_{2m}$ ($m\ge 2$)  be the set whose elements  are  all cycles $h_{2m}$ in $G$ of length $2m$ which are properly half-monochromatic by the coloring $\mathcal C$.
Finally, let $X= K\cup(\bigcup_{m\ge 2}B_{2m})\cup (\bigcup_{m\ge 2}H_{2m})$.
Again the elements of $X$ are (some of) the subsets of $E$.

Let us now recolor each edge $e\in E$ using a new color randomly and independently with probability
${c\over \D}$ (where $c\le 1$ is a constant to be determined later). Call $\mathcal C'$ this new coloring which by construction uses $\D+2$ colors.
We need to  show  that with positive probability the coloring  $\mathcal C'$ is such that

\begin{itemize}
\item[A.] No pair of adjacent edges has the same color

\item[B.] There is no properly bichromatic cycle
 \end{itemize}

\\We use condition \equ(eq:r4) of Theorem  \ref{LLLn}. So, once $G$ has been recolored by the coloring $\mathcal C'$, consider the following bad events.

 \begin{itemize}
\item[I.] For each pair of adjacent edges $\{e,e'\}\in K$, let ${ A}_{\{e,e'\}}$ be the
event that $e$ and $e'$ have the same color.

\item[II.] For each properly bichromatic cycle $b_{2k}\in B_{2k}$ of length $2k$, $k\ge 2$, in $G$ with respect to the coloring $\mathcal C$, let
$A_{b_{2k}}$ be the event that either no edge is recolored with the
new color or one half is recolored and the other half stays unchanged.

\item[III.] For each properly half-monochromatic cycle $h_{2m}\in H_{2m}$ of length $2m$, $m\ge 2$,  with respect to $\mathcal C$,
let ${A}_{h_{2m}}$ be the event that  $h_{2m}$  becomes
properly bichromatic by recoloring
the non monochromatic  half of its edges with the new color and leaving the monochromatic part unchanged.

\end{itemize}

Clearly, if none of the bad events I, II, or III occurs, properties A and B are satisfied. It is straightforward to see that the probabilities of
the events I, II, and III are as follows.

For each pair of adjacent edges $\{e,e'\}$
$$
{\rm Prob}({ A}_{\{e,e'\}})= {c^2\over \D^2}
$$
For $k\ge 2$ and each properly bichromatic cycle $b_{2k}$ of length $2k$,
$$
{\rm Prob}(A_{b_{2k}})=\left(1-{c\over \D}\right)^{2k}+
2\left(1-{c\over \D}\right)^{k}\left({c\over \D}\right)^{k}\le
$$
$$
\le {1\over \left(1+{c\over \D}\right)^{2k}}~~~~~~~~~~~~~~~~\Eq(A)
$$
For $k\ge 2$ and each properly half-monochromatic cycle $h_{2k}$ of length $2k$,
$$
{\rm Prob}(A_{h_{2k}})= {c^k\over \D^k}\left(1-{c\over \D}\right)^{k}
$$
To prove \equ(A) let  $w=c/\D$. Then  \equ(A) becomes
$$
\left(1-w\right)^{2k}+ 2w^{k}\left(1-w\right)^{k}\le {1\over \left(1+w\right)^{2k}}
$$
i.e., multiplying  both side of the inequality by $(1+w)^{2k}$
$$
\left(1-w^2\right)^{2k}+ 2\left(1-w^2\right)^{k}\left(w(1+w)\right)^{k}\le {1}
$$
which is true for all $k\ge 2$ , if $w\le {1/2}$, which is indeed the case since $w=c/\D\le 1/3$.

Now, as before,  a bad event $A_x$ of the collection $\{A_z\}_{z\in X}$ is independent of all other events $A_{x'}$
such that  $x\cap x'=\emptyset$. So the graph $H=(X,F)$  with vertex
set $X$ and edge set $F$ such that
the pair $\{x,x'\}\in F$ if and only if  $x\cap x'\neq\emptyset$ is a dependency graph for the collection of events $\{A_x\}_{x\in X}$.
Moreover, as shown in \cite{ASZ} (see there Lemma 7),  in a properly edge-colored graph $G$

\begin{itemize}
\item
each edge $e$ is contained in at most $2\D$ pairs $\{f,f'\}$ of incident edges
\item
each edge $e$ is contained in at most $\D$ properly bichromatic cycles of $G$

\item
each edge $e$ is contained in at most $2\D^{k-1}$ half-monochromatic cycles of length $2k$
\end{itemize}
Hence

\begin{itemize}
\item[{[a]}]
For each  vertex $x=\{e_1,e_2\}\in K$ of $H$,
$\G_H^*(x)$ is the union of $2$ sets $\G_i^*(x)$ ($i=1,2$)
such that
$$
|\G_i^*(x)|\le 2{\D} +
\Delta + \sum_{k\ge \lceil g/2\rceil }2\D^{k-1}
$$
and every element of $\G_i^*(x)$ contains $e_i$.
Hence  the subgraph of $H$  induced by $\G_i^*(x)$ is a clique for $i=1,2$.
\item[{[b]}]
For $k\ge 2$ and for each  vertex $y=\{e_1,\dots,e_{2k}\}\in B_{2k}\cup H_{2k}$ of $H$, we have that
$\G_H^*(y)$ is the union of $m$ sets $\G_1^*(y),\dots,\G_{2k}^*(y)$ such that,
for $j=1,\dots,  2k$,
$$
|\G_j^*(y)|\le 2{\D} +
\Delta + \sum_{k\ge \lceil g/2\rceil }2\D^{k-1}
$$
and every element of $\G_{j}^*(y)$ contains $e_j$ so that the subgraph of $H$  induced by $\G_j^*(y)$ is a clique
for all $j=1,2,\dots, 2k$.
\end{itemize}
Let us now choose nonnegative  numbers $\{\m_z\}_{z\in X}$as follows. For any $x\in K$, put $\m_{x}=\m_1$; for
each $y\in B_{2m}$, put $\m_y=\m_2$;  for any $z\in H_{2m}$ put
$\m_{z}=\m_{2m}$. Then using once again inequality \equ(key) one gets, under conditions [a] and [b], that, for any $x\in K$ and any $y\in B_{2k}\cup H_{2k}$,

$$
\ffp_{x}(\bm \mu)\le(1+ 2\D\m_{1} + \D\m_2+\sum_{s\ge \lceil g/2\rceil }2\D^{s-1}\m_{2s})^{2}
$$
and
$$
\ffp_{y}(\bm \mu)\le(1+ 2\D\m_{1} + \D\m_2+\sum_{s\ge \lceil g/2\rceil }2\D^{s-1}\m_{2s})^{2k}
$$

Hence, we have that the condition of Theorem \ref{LLLn} is satisfied
if there are positive numbers $\m_1,\m_2,\{\m_{2k}\}_{k\ge  \lceil g/2\rceil}$ such that the following inequalities are satisfied
$$
{c^2\over \D^2}\le {\m_{1}\over (1+ 2\D\m_1 + \D\m_2+ \sum_{s\ge \lceil g/2\rceil}2\D^{s-1}\m_{2s})^2}
$$

$$
{1\over \left(1+{c\over \D}\right)^{2k}}\le
{\m_{2}\over (1+ 2\D\m_1 + \D\m_2+ \sum_{s\ge \lceil g/2\rceil}2\D^{s-1}\m_{2s})^{2k}}
$$

$$
{c^k\over \D^k}\left(1-{c\over \D}\right)^{k} \le {\m_{2k}
\over (1+ 2\D\m_1 + \D\m_2+ \sum_{s\ge \lceil g/2\rceil}2\D^{s-1}\m_{2s})^{2k}}
$$
Put now $\m_1=\m_2={\a^2\over \D^2}$ and $\m_{2k}={\a^k\over \D^k}$. Then the inequalities above become
$$
c ~\le~ {\a\over 1+ {R_g(\a)\over \D}}\Eq(ine1)
$$

$$
{1\over \left(1+{c\over \D}\right)^k}~\le~
{{\a\over\D}\over \left[1+ {R_g(\a)\over \D}\right]^k}\Eq(ine2)
$$

$$
c \left(1-{c\over \D}\right)~ \le ~{\a\over \left[1+{R_g(\a)\over \D}\right]^{2}}\Eq(ine3)
$$

where we have put
$$
R_g(\a)= 3\a^2+ {2\a^{\lceil g/2\rceil}\over 1-\a}\Eq(rga)
$$
Now note that \equ(ine2) can be satisfied for all $k$ greater than some fixed $k_0$ only if
$$
c> R_g(\a)\Eq(strictu)
$$
On the other hand, if \equ(strictu) holds, then  inequality \equ(ine3) is  satisfied if \equ(ine1) is satisfied. Indeed inequality \equ(ine3)  can be rewritten as
$$
{c} \le {\a\over \left(1+ {R_g(\a)\over \D}\right)^{2} \left(1-{c\over \D}\right)}=  {\a\over \left(1+ {R_g(\a)\over \D}\right)}{1\over  \left(1-{c\over \D}\right) \left(1+ {R_g(\a)\over \D}\right)}
$$
and $\left[\left(1-{c\over \D}\right) \left(1+ {R_g(\a)\over \D}\right)\right]^{-1}>1$ due to \equ(strictu).

\\Now let us suppose that $g\ge 80$  so that
$$
R_{g}(\a)\le 3\a^2+ {2\a^{40}\over 1-\a}
$$
and let find $\a_0$ such that
$$
f(\a_0)={\a\over \left(1+ {R_{80}(\a)\over \D}\right)} -R_{80}(\a)
$$
is maximum. A simple calculation show that the maximum occurs around 0$.155$ so let us choose $\a=\a_0= 0.155$ and $f(\a_0)= 0.07928$. Hence we can choose
$$
c~= ~c_0~={\a_0\over \left(1+ {R_{80}(\a_0)\over \D}\right)}
$$
Now
solve
$$
{1\over \left(1+{c_0\over \D}\right)^k}\le {{\a_0\over \D}\over \left(1+ {R_{80}(\a_0)\over \D}\right)^k}
$$
i.e.
$$
k~\log\left\{1+{c_0\over \D}\over 1+ {R_{80}(\a_0)\over \D}\right\}\ge \log\D + \log(1/\a_0)
$$
Observe now that, for $a>b>0$, $\log\{{1+a\over 1+ b}\}= \log\{1+{a-b\over 1+ b}\}$ and,  for  $w\ge 0$,  $\log(1+w)\ge w(1-w/2)$. So,
recalling that $R_{80}(\a_0)=0.0721$, $c_0- R_{80}(\a_0)= 0.07928$, and  also using that  $\D\ge 3$, we get
after some easy computations
$$
k\ge  12.92\D\log\D\Big(1+ {2\over \log\D}\Big) \Big(1 + {2\over \D}\Big)
$$
which holds for all integers $k\ge  \lceil g/2\rceil$ as soon as $g\ge  25.84\D\log\D(1+ {2\over \log\D}) (1 + {2\over \D})$  and since
$(1+ {2\over \log\D}) (1 + {2\over \D}) \le 1+ {4.1\over \log \D}$ for all $\D\ge 3$, Theorem 3 item (c) follows. $\Box$
\v\v


\subsection{Proof item (d): acyclic chromatic number of $G$}

\def\N{{\mathbb N}}\def\la{\langle}\def\ra{\rangle}
We follow  \cite{AMR} using Theorem \ref{LLLn} instead of the Lov\'asz Local Lemma.
Let  $\mathcal C$ be a  vertex-coloring of $G=(V,E)$  such that in  each vertex
the color is chosen at random independently and uniformly among
$N\ge c\D^{4/3}$ colors ($c$ being a positive constant to be determined later).
In the following a pair of non-adjacent vertices $u,v$ of $G$ will be  called a {\it special pair} if
$u$ and $v$ have more than $\D^{2/3}$ common neighbors and will be denoted by $\la u,v\ra$.

Let $P_{4}$ be the set whose elements  are  set of vertices  $\{v_0,v_1,v_2,v_3,v_4\}$ forming paths of length four  in $G$. Let $C_4$
be the set whose elements  are  sets of vertices  $\{v_1,v_2,v_3,v_4\}$ forming 4-cycles  in $G$. Let $S$  be the set whose elements  are  sets of vertices  $\la v,v'\ra$ forming special pairs  in $G$.
Finally, let $X= E\cup P_4\cup C_4\cup S$ (here, of course, $E$ is the set of edges of $G$). Observe that
now the elements of $X$ are (some of) the subsets of  the vertex set $V$ of $G$.

Consider the following unfavorable events.

 \begin{itemize}
\item[I.] For each pair of adjacent vertices $\{u,v\}\in E$ of $G$, let $A_{\{u,v\}}$ be the
event that $u$ and $v$ have the same color.

\item[II.] For each path  $p_4=v_0v_1v_2v_3v_4\in P_4$ of $G$, let
$A_{p_4}$ be the event that vertices  $v_0, v_2,v_4$ have the same color and  vertices  $v_1,v_3$ have the same color.

\item[III] For each induced 4-cycle $c_4=v_1v_2v_3v_4\in C_4$ of $G$, in which neither $v_1; v_3$
nor $v_2; v_4$ is a special pair, let $A_{c_4}$be the event that  $v_1; v_3$ have the same color and
$v_2; v_4$ have the same color.

\item[IV] For each special pair of vertices $\la u,v\ra\in S$ of $G$ let $A_{\la u,v\ra}$ be the
event that $u$ and $v$ receive the same color.

\end{itemize}

Alon, Mc Diarmid and  Reed   have shown in \cite{AMR}
that if none of the event I, II,  III or IV occurs then the graph is properly colored without bichromatic cycles (see in \cite{AMR}, proof of proposition 2.2).
We now use Theorem \ref{LLLn} to show that with positive probability  none of the events
occurs.

We first observe that  the probability
of an  event of type I, II, III and IV respectively  are
$$
{\rm Prob}(A_{\{u,v\}})= {1\over N},~~~~
{\rm Prob} (A_{p_4})= {1\over N^3}
~~~~~~
{\rm Prob}(A_{c_4})= {1\over N^2}
~~~~~~
{\rm Prob} (A_{\la u,v\ra})= {1\over N}
$$

Secondly, we note that, for $x,x'\in X$ an event $A_x$ (where $x\in X$ can be a pair of adjacent vertices, a path of length four,
a cycle of length three or a special pair)
is independent of all events $A_{x'}$ such that $x\cap x'=\emptyset$. So the graph $H=(X,F)$  with vertex
set $X$ and edge set $F$, such that
the pair $\{x,x'\}\in F$ if and only if  $x\cap x'\neq\emptyset$, is a dependency graph for the collection of events $\{A_x\}_{x\in X}$.

Finally, following \cite{AMR}, (see there the proof of Lemma 2.4) we have that

\begin{itemize}
\item
a vertex $v\in V$ belongs to at most $\D$ edges of $G$.
\item
a vertex $v\in V$ belongs to at most 
$ {5\over 2}\D^4$  paths of length 4 in $G$.

\item
The number of induced 4-cycles in $G$ containing $v$ in which no opposite
pair of vertices is a special pair is at most ${1\over 2}\D^{8/3}$.
\item
 The number of special pairs of vertices containing a given vertex $v$ is at most $\D^{4/3}$.
\end{itemize}

Hence

\begin{itemize}
\item[{[a]}]
For each  vertex $x\in E\cup S$ of $H$ (i.e. either $x=\{v_1,v_2\}$ or $x=\la v_1,v_2\ra$),
$\G_H^*(x)$ is the union of $2$ sets $\G_i^*(x)$ ($i=1,2$)
such that
$$
|\G_i^*(x)|\le  \D + {5\over 2}\D^4+ {1\over 2}\D^{8/3}+\D^{4/3}
$$
and every element of $\G_i^*(x)$ contains $v_i$.
Hence  the subgraph of $H$  induced by $\G_i^*(x)$ is a clique for $i=1,2$.
\item[{[b]}]
For  each  vertex $y= \{v_0,v_1,v_2,v_3,v_4\}\in P_4$ of $H$, we have that
$\G_H^*(y)$ is the union of $5$ sets $\G_1^*(y),\dots,\G_{5}^*(y)$ such that,
for $j=1,\dots,  5$,
$$
|\G_j^*(y)|\le \D +{5\over 2}\D^4+ {1\over 2}\D^{8/3}+\D^{4/3}
$$
and every element of $\G_{j}^*(y)$ contains $v_j$ so that the subgraph of $H$  induced by $\G_j^*(y)$ is a clique
for all $j=1,2,\dots, 5$.
\item[{[c]}]
For  for each  vertex $z= \{v_1,v_2,v_3,v_4\}\in C_4$ of $H$, we have that
$\G_H^*(z)$ is the union of $4$ sets $\G_1^*(z),\dots,\G_{4}^*(z)$ such that,
for $j=1,\dots,  4$,
$$
|\G_j^*(z)|\le \D + {5\over 2}\D^4+ {1\over 2}\D^{8/3}+\D^{4/3}
$$
and every element of $\G_{j}^*(z)$ contains $v_j$ so that the subgraph of $H$  induced by $\G_j^*(z)$ is a clique
for all $j=1,2,3,4$.
\end{itemize}
Let us now choose nonnegative  numbers $\{\m_u\}_{u\in X}$ as follows. For any $x\in E$, put $\m_{x}=\m_1$; for
each $y\in P_4$, put $\m_y=\m_2$; for any $z\in C_4$, put
$\m_{z}=\m_{3}$;  for any $w\in S$, put $\m_{z}=\m_{4}$.
Then using  \equ(key) one gets, under conditions [a]-[c], that, for any $x\in E\cup S$, any $y\in P_4$, and any $z\in C_4$

$$
\ffp_{x}(\bm \mu)~\le~(1+ \D\m_1 + {5\over 2}\D^4\m_2+ {1\over 2}\D^{8/3}\m_3+\D^{4/3}\m_4)^2
$$
and
$$
\ffp_{y}(\bm \mu)~\le~(1+ \D\m_1 + {5\over 2}\D^4\m_2+ {1\over 2}\D^{8/3}\m_3+\D^{4/3}\m_4)^5
$$
and
$$
\ffp_{z}(\bm \mu)~\le~(1+ \D\m_1 + {5\over 2}\D^4\m_2+ {1\over 2}\D^{8/3}\m_3+\D^{4/3}\m_4)^4
$$

Hence Theorem \ref{LLLn} holds if
there are positive numbers $\m_1$, $\m_2$, $\m_3$, and $\m_4$ such that
the following inequalities are simultaneously satisfied:

$$
{1\over N}\le {\m_{1}\over (1+ \D\m_1 + {5\over 2}\D^4\m_2+ {1\over 2}\D^{8/3}\m_3+\D^{4/3}\m_4)^2}
$$

$$
{1\over N^3}\le {\m_{2}\over (1+ \D\m_1 + {5\over 2}\D^4\m_2+ {1\over 2}\D^{8/3}\m_3+\D^{4/3}\m_4)^5}
$$

$$
{1\over N^2}\le {\m_{3}\over (1+ \D\m_1 + {5\over 2}\D^4\m_2+ {1\over 2}\D^{8/3}\m_3+\D^{4/3}\m_4)^4}
$$

$$
{1\over N}\le {\m_{4}\over (1+ \D\m_1 + {5\over 2}\D^4\m_2+ {1\over 2}\D^{8/3}\m_3+\D^{4/3}\m_4)^2}
$$
Taking $\m_1=\m_4=\m$,  $\m_2=\m^3$ and $\m_3=\m^2$ these inequalities are satisfied if, for some $\m>0$

$$
{1\over N}\le {\m\over (1+ (\D+ \D^{4/3}) \m + {1\over 2}\D^{8/3}\m^2 + {5\over 2}\D^4\m^3)^2}
$$
Now choose $\m=\a/\D^{4/3}$. Then, recalling that $N\ge c\D^{4/3}$, inequality above is satisfied if
$$
{1\over c}\le {\a\over(1+ (1+ \D^{-1/3}) \a + {1\over 2}\a^2 + {5\over 2}\a^3)^2}
$$
i.e.
$$
c\ge {1\over \a}(1+  \a + {1\over 2}\a^2 + {5\over 2}\a^3)^2 + \left[ {\a\over \D^{2/3}}+ {2\over \D^{1/3}}(1+  \a + {1\over 2}\a^2 + {5\over 2}\a^3)\right]\Eq(avb)
$$
and taking  $\a=0.34$  it is easy  to check that the right hand side of \equ(avb) is less  than $6.583+ 3.3/ \D^{1/3}$ for all $\D\ge 3$.
So we get $c\ge 6.583+ 3.3/ \D^{1/3}$ and hence $a(G)\le \lceil 6.583 \Delta^{4/3}+ 3.3\D\rceil$. $\Box$

\def\N{{\mathbb N}}

\subsection{Proof of item (e): star  chromatic number of $G$}

We follow \cite{FRR}, but  we use Theorem \ref{LLLn} in place of the Lov\'asz Local Lemma.
Let   $\mathcal C$ be a vertex-coloring of $G=(V,E)$ using $N\ge c\D^{3/2}$ colors ($c$ being a positive constant to be determined later)
such that in  each vertex the color
is chosen at random independently and uniformly among the  set of
$N$ colors.

Let $P_{3}$ be the set whose elements  are  set of vertices  $\{v_1,v_2,v_3,v_4\}$ forming paths of length three  in $G$ and
 let $X= E\cup P_3$. Observe that, as in subsection 3.4, the elements of $X$ are (some of) the subsets of  the vertex set $V$ of $G$.

Consider the following unfavorable events.

 \begin{itemize}
\item[I.] For each pair of adjacent vertices $\{u,v\}\in E$ of $G$, let $A_{\{u,v\}}$ be the
event that $u$ and $v$ have the same color.

\item[II.] For each path of length three $p_3=v_1v_2v_3v_4\in P_3$ in $G$, let
$A_{p_3}$ be the event that vertices $v_1$, $v_3$ have  the same color and  vertices $v_2$ and $v_4$ have the same color.

\end{itemize}
Clearly, by definition, if none of the events above occurs then $\mathcal C$ is a star coloring.

We first observe that  the probability
of an  event of type I, II,  respectively  are
$$
{\rm Prob}(A_{\{u,v\}})= {1\over N},~~~~
{\rm Prob}(A_{p_3})= {1\over N^2}
$$
Then, as in Subsection 3.4, we observe that, for $x\in X$, the event $A_x$ (where now $x\in X$ can be a pair of adjacent vertices or  a path of length four)
is independent of all events $A_{x'}$, with $x'\in X$, such that $x\cap x'=\emptyset$. So the graph $H=(X,F)$  with vertex
set $X$ and edge set $F$ such that
the pair $\{x,x'\}\in F$ if and only if  $x\cap x'\neq\emptyset$ is a dependency graph for the collection of events $\{A_x\}_{x\in X}$.

Finally, as observed in \cite{FRR} (see there Observation 8.1)  we have that

\begin{itemize}
\item
a vertex $v\in V$ belongs to at most $\D$ edges of $G$.
\item
a vertex $v\in V$ belongs to at most $2\D(\D-1)^2\le 2\D^3$  paths of length 3 in $G$.
\end{itemize}
Hence
\begin{itemize}
\item[{[a]}]
For each vertex  $x=\{v_1,v_2\}\in E$ of $H$,
$\G_H^*(x)$ is the union of $2$ sets $\G_i^*(x)$ ($i=1,2$)
such that
$$
|\G_i^*(x)|\le  \D + 2\D^3
$$
and every element of $\G_i^*(x)$ contains $v_i$.
Hence  the subgraph of $H$  induced by $\G_i^*(x)$ is a clique for $i=1,2$.
\item[{[b]}]
For  for each  vertex $y= \{v_1,v_2,v_3,v_4\}\in P_3$ of $H$, we have that
$\G_H^*(y)$ is the union of $4$ sets $\G_1^*(y),\dots,\G_{4}^*(y)$ such that,
for $j=1,\dots,  4$,
$$
|\G_j^*(y)|\le \D + 2\D^3
$$
and every element of $\G_{j}^*(y)$ contains $v_j$ so that the subgraph of $H$  induced by $\G_j^*(y)$ is a clique
for all $j=1,2,\dots, 4$.
\end{itemize}
Let us now choose nonnegative  numbers $\{\m_z\}_{z\in X}$as follows. For any $x\in E$, put $\m_{x}=\m_1$; for
each $y\in P_3$, put $\m_y=\m_2$.
Then, by    conditions [a] and [b], using \equ(key) we get that, for any $x\in E$, and any $y\in P_3$

$$
\ffp_{x}(\bm \mu)~\le~(1+ \D\m_1 + 3\D^3\m_2)^2
$$
and
$$
\ffp_{y}(\bm \mu)~\le~(1+ \D\m_1 + 3\D^3\m_2)^4
$$
Hence, analogously to the previous sections, Theorem \ref{LLLn} holds  if we can find
nonnegative  numbers $\m_{1}$, $\m_{2}$
such that the following inequalities are simultaneously satisfied:

$$
{1\over N}\le {\m_{1}\over (1+ \D\m_1 + 2\D^3\m_2)^2}
$$

$$
{1\over N^2}\le {\m_{2}\over (1+ \D\m_1 + 2\D^3\m_2)^4}
$$
Now take  $\m_2=\m_1^2$ and $\m_1=\a/\D^{3/2}$. Then, recalling that $N\ge c\D^{3/2}$, these inequalities are
satisfied if, for some $\a>0$

$$
{1\over c}\le {\a\over(1+ \D^{-1/2}\a + 2\a^2)^2}\Eq(above)
$$
Maximizing the right hand side of \equ(above) with respect to $\a$  we get that the maximum is reached at
$$
\a=\a_0= {1\over \sqrt{6}\left(~\sqrt{1+{1\over 24\D}}+ \sqrt{1\over 24\D}~\right)} $$
Now observing that $\a_0\le 1/\sqrt{6}$ for all $\D\ge 1$, we get   that the inequality \equ(above) is satisfied for all $\D\ge 3$  as soon as
$$
c\ge  \sqrt{6}\left[~\sqrt{1+{1\over 24\D}}+ \sqrt{1\over 24\D}~\right]\left[{4\over 3}+{1\over \sqrt{6\D}}\right]^2\Eq(abv)
$$
It is now  easy to check that the left hand side of \equ(abv) is less than ${16\over 9}\sqrt{6} + {1.5\over\sqrt{\D}}$ for all $\D\ge 3$. So we get
$c\ge {16\over 9}\sqrt{6} + {1.5\over\sqrt{\D}}$ and hence $\chi_s(G)\le \lceil4.34\Delta^{3/2}+ 1.5\D\rceil$.
$\Box$

\subsection{Proof of item (f): $\b$-frugal chromatic number}
We follow  \cite{HMR}, but we  use Theorem \ref{LLLn} instead of the Lov\'asz Local Lemma.
Let now  $\mathcal C$ be a vertex-coloring of $G=(V,E)$ using $c$ colors ($c$ being a positive constant to be determined later)
such that in  each vertex the color
is chosen at random independently and uniformly among the  set of
$c$ colors.

We may assume $\b\ge 2$ since in the case $\b=1$
the $1$-frugal chromatic number of $G$,
$\chi^1(G)$, is just the vertex chromatic number of the graph obtained from $G$ by adding an edge between
any two vertices at distance 2 in $G$,
which has maximum degree at most $\D^2$ and hence, by Vizing
$\chi^1(G)\le \D^2+1$. Given $v\in V$, let $S^v_{\b}$ be the set whose elements  are  sets of vertices  $\{v_1,\dots ,v_{\b+1}\}$ such that
$\{v_1,\dots ,v_{\b+1}\}\subset  \G_G(v)$. Let $S_\b=\bigcup_{v\in V}S^v_{\b}$.
Let $X= E\cup S_{\b}$. As in subsections 3.4 and 3.5, the elements of $X$ are (some of) the subsets of  the vertex set $V$ of $G$.

 Consider the following unfavorable events.
\begin{itemize}
\item[I.]
For $\{u,v\}\in E$, let
$A_{\{u,v\}}$ be the event that $u$ and $v$ receive the same color.
\item[II.]
For
$s_\b=\{v_1,\dots, v_{\b+1}\}\in S_{\b}$,
let $A_{s_{\b}}$ be the event that all vertices in $s_\b$ receive the same color.
\end{itemize}
If none of the events above occurs then, by definition, $\mathcal C$ is a $\b$ frugal coloring. We have clearly
$$
{\rm Prob} (A_{\{u,v\}}) = {1\over c}
~~~~~~~~{\rm and}~~~~~~~~
{\rm Prob} (A_{h_{\b}}) = {1\over c^\b}
$$
As in subsections 3.4 and 3.5, the graph $H=(X,F)$  with vertex
set $X$ and edge set $F$ such that
the pair $\{x,x'\}\in F$ if and only  $x\cap x'\neq\emptyset$ is a dependency graph for the events $\{A_x\}_{x\in X}$.

Moreover
\begin{itemize}
\item
a vertex $v\in V$ belongs to at most $\D$ edges of $G$.
\item
a vertex $v\in V$ belongs to at most $\D {\D\choose\b}\le {\D^{1+\b}/ \b!}$
sets of type $h_\b=\{v_1,\dots, v_{\b+1}\}$.
\end{itemize}

Hence
\begin{itemize}
\item[{[a]}]
For each  vertex $x=\{v_1,v_2\}\in E$ of $H$,
$\G_H^*(x)$ is the union of two sets $\G_i^*(x)$ ($i=1,2$)
such that
$$
|\G_i^*(x)|\le  \D + {\D^{1+\b}/ \b!}
$$
and every element of $\G_i^*(x)$ contains $v_i$.
Hence  the subgraph of $H$  induced by $\G_i^*(x)$ is a clique for $i=1,2$.
\item[{[b]}]
For  for each  vertex $y=\{v_1,\dots,v_{\eta+1}\}\in S_{\b}$ of $H$, we have that
$\G_H^*(y)$ is the union of $\b+1$ sets $\G_1^*(y),\dots,\G_{\b+1}^*(y)$ such that,
for $j=1,\dots,  \b+1$,
$$
|\G_j^*(y)|\le \D +{\D^{1+\b}/ \b!}
$$
and every element of $\G_{j}^*(y)$ contains $v_j$ so that the subgraph of $H$  induced by $\G_j^*(y)$ is a clique
for all $j=1,2,\dots, \b+1$.
\end{itemize}

Let us now choose nonnegative  numbers $\{\m_z\}_{z\in X}$ as follows. For any $x\in E$, put $\m_{x}=\m_1$;   for
each $y\in S_\b$, put $\m_y=\m_2$.
Then, by    conditions [a] and [b], using \equ(key) we get that, for any $x\in E$ and any $y\in S_\b$

$$
\ffp_{x}(\bm \mu)~\le~(1+ \D\m_1 +{\D^{1+\b}\over \b!}\m_2)^{2}
$$
and
$$
\ffp_{y}(\bm \mu)~\le~(1+ \D\m_1 +{\D^{1+\b}\over \b!}\m_2)^{\b+1}
$$
Hence
Theorem \ref{LLLn} holds if
$$
{1\over c}\le {\m_1\over (1+ \D\m_1 + {\D^{1+\b}\over \b!}\m_2)^2} \Eq(81)
$$
and
$$
{1\over c^\b}\le  {\m_2\over (1+ \D\m_1 +{\D^{1+\b}\over \b!}\m_2)^{\b+1}} \Eq(82)
$$
Put $\m_1= \m$, $\m_2=\b!\m^{1+\b}$ and $\a=\D\m\in (0,+\infty)$, then inequalities \equ(81) and \equ(82) become

$$
{1\over c}\le {1\over \D} {\a\over (1+ \a + \a^{1+\b})^2} \Eq(81a)
$$
and
$$
{1\over c}\le {(\b!)^{1/\b}\over \D^{1+{1\over\b}}}  \left[{\a\over (1+ \a + \a^{1+\b})}\right]^{1+{1\over \b}}\Eq(82b)
$$
which are satisfied if
$$
c\ge \max\Big\{k_1(\b)\D\,,\, \,k_2(\b){\D^{1+{1\over \b}}\over (\b!)^{1/\b}}\Big\}\Eq(upchi)
$$
where
$$
k_1(\b)= \min_{\a>0} {(1+ \a + \a^{1+\b})^2\over \a}
$$
$$
k_2(\b)= \left[\min_{\a>0}  {1+ \a + \a^{1+\b}\over \a}\right]^{1+{1\over \b}}
$$
An easy computation shows that $k_1(\b)$ and $k_2(\b)$ are both decreasing functions of $\b$
and,  for $\b\ge 2$, we have that
$k_1(\b)\in[4, 5.27]$
and
$k_2(\b)\in[2,4.92]$. Therefore we get that $c\ge \max\{k_1(\b)\D\,,\, \,k_2(\b){\D^{1+{1/ \b}}/ (\b!)^{1/\b}}\}$ and hence
$\chi^\b(G)\le \lceil\max\{k_1(\b)\D, \; k_2(\b){\D^{1+{1/ \b}}/ (\b!)^{1/\b}}\}\rceil$.
$\Box$

\vskip.2cm
\\{\bf Remark}. Note that \equ(upchi) is valid for all $\D\ge 3$. Moreover the upper bound \equ(upchi) of $\chi^\b(G)$ is a decreasing function of $\b$.

\section*{Acknowledgments}

This work  has been partially supported by
Conselho Nacional de Desenvolvimento Cient\'{i}fico e Tecnol\'ogico (CNPq),
 CAPES (Coordena\c{c}\~ao de Aperfei\c{c}oamento de Pessoal de
N\'{\i}vel Superior, Brasil), FAPEMIG (Funda{c}\~ao de Amparo \`a Pesquisa do Estado de Minas Gerais)
and Universit\`a di Roma ``Tor Vergata".

\end{document}